\documentclass{amsart}
\usepackage{amsfonts}
\usepackage{amscd}
\usepackage{amssymb}
\newtheorem{theorem}{Theorem}
\theoremstyle{plain}

\newtheorem{corollary}[theorem]{Corollary}

\newtheorem{lemma}[theorem]{Lemma}

\newtheorem{proposition}[theorem]{Proposition}

\theoremstyle{definition}

\newtheorem{remark}[theorem]{Remark}

\newcommand{\R}{{\mathbb R}}

\newcommand{\C}{{\mathbb C}}
\newcommand{\N}{{\mathbb N}}

\renewcommand{\H}{{\mathcal H}}

\renewcommand{\S}{{\mathcal S}}
\renewcommand{\Re}{{\mathop{\mathrm{Re}}\,}}
\renewcommand{\Im}{{\mathop{\mathrm{Im}}\,}}

\newcommand{\supp}{{\operatorname {supp}\,}}
\normalbaselines
\begin{document}

\title[The Paley--Wiener theorem]
{{Elementary proofs of Paley--Wiener theorems for the Dunkl transform on the real
line}}
\author{Nils Byrial Andersen}
\address{Dipartimento di Matematico,
Politecnico di Torino,
Corso Duca degli Abruzzi, 24,
10129 Torino, Italy}
\email{byrial@calvino.polito.it}
\thanks{The first author is supported by a research grant from the
{\it{European Commission}} IHP Network: 2002--2006
{\it{Harmonic Analysis and Related Problems}}
(Contract Number: HPRN-CT-2001-00273 - HARP)}
\author{Marcel de Jeu}
\address{Mathematical Institute,
Leiden University,
P.O. Box 9512,
2300 RA Leiden,
The Netherlands}
\email{mdejeu@math.leidenuniv.nl}

\begin{abstract}
We give an elementary proof of the Paley--Wiener theorem for smooth functions
for the Dunkl transforms on the real line, establish a similar theorem for
$L^2$-functions
and prove identities in the spirit of Bang for $L^p$-functions. The proofs seem to
be new also in the special case of the Fourier transform.
\end{abstract}

\maketitle

\section{Introduction and overview}

The Paley--Wiener theorem for the Dunkl transform $D_k$ with multiplicity $k$
(where $\Re k\ge 0$) on the real line states that a smooth function $f$ has
support in the bounded interval $[-R,R]$ if, and only if, its transform $D_k f$
is an entire function which satisfies the usual growth estimates as they are
required in the (special) case of the Fourier transform. Various proofs of this
result are known, all of which use explicit formulas available in this
one-dimensional setting (see Remark~\ref{PWremark} for more details). In this
paper, however, we present an alternative proof which does not use such
explicit expressions, being based almost solely on the formal properties of the
transform. Along the same lines, we also obtain a Paley--Wiener theorem for
$L^2$-functions for $k\ge 0$. The case $k=0$ specializes to the Fourier
transform, and to our knowledge the proofs of both Paley--Wiener theorems are
new even in this case.

In addition, we establish two identities in the spirit of Bang \cite[Theorem
1]{Ba1}. These results could be called real Paley--Wiener theorems (although
terminology is not yet well-established), since they relate certain growth
rates of a function \emph{on the real line} to the support of its transform.
The approach at this point is inspired by similar techniques in \cite{An1, An2,
An3}. Our results in this direction partially overlap with \cite{CT, COT}, but
the new proofs are considerably simpler, as they are again based almost solely
on the formal properties of the transform. We will comment on this in more
detail later on, as these results will have been established. For $k=0$, one
retrieves Bang's result; we feel that the present method of proof, which, e.g.,
does not use the Paley--Wiener theorem for smooth functions, but rather implies
it, is then more direct than that in \cite{Ba1}.

The rather unspecific and formal structure of the proofs suggests that the
methods can perhaps be put to good use for other integral transforms with a
symmetric kernel, both for the Paley--Wiener theorems and the equalities in the
spirit of Bang (cf.\ Remark~\ref{PWremark}). The structure of the proof is also
such that, if certain combinatorial problems can be surmounted, a proof of the
Paley--Wiener theorem for the Dunkl transform for invariant balanced compact
convex sets in arbitrary dimension might be possible. This would be further
evidence for the validity of this theorem for invariant compact convex sets
(cf.~\cite[Conjecture 4.1]{dJ2}), but at the time of writing this higher-dimensional result has not been established.

This paper is organized as follows.
In Section~\ref{basics} the necessary notations and previous results are given.
Section~\ref{PWfunctions} contains the Paley--Wiener theorem for smooth functions
and can---perhaps---serve as a model for a proof of such a theorem in other
contexts. The rest of the paper is independent of this section. Section~\ref{Lp} is
concerned with the real Paley--Wiener theorem in the $L^p$-case and the $L^2$-case
is settled in Section~\ref{L2}.

\section{Dunkl operators and the Dunkl transform on $\R$}\label{basics}

The Dunkl operators and the Dunkl transform were introduced for arbitrary root
systems by Dunkl \cite{Du1, Du2, Du3}. In this section we recall some basic
properties for the one-dimensional case of $A_1$, referring to \cite[Sections~1
and 2]{roesler} for a more comprehensive overview and to \cite{Du1, Du2, Du3,
dJ1, besseletc, roeslerpositivity} for details. We suppress the various
explicit formulas which are known in this one-dimensional context (as these are
not necessary for the proofs), thus emphasizing the basic structure of the
problem which might lead to generalizations to the case of arbitrary root
systems.

Let $k \in \C$, and consider the Dunkl operator $T_k$
\begin{equation*}
T_k f(x) = f'(x) +k  \frac {f(x) - f(-x)}x\qquad (f\in C^\infty(\R),
\,x\in \R).
\end{equation*}
Rewriting this as
\begin{equation}
T_k f(x) = f'(x) +  k
\int _{-1} ^1 f'(tx)\, dt,
\label{intDunkl}
\end{equation}
it follows that $T_k$ maps $C^\infty (\R)$, $C_c^\infty (\R)$ and the Schwartz space
$\S (\R)$ into themselves.

If $\Re k\geq 0$, as we will assume for the remainder of this section, then, for each
$\lambda\in\C$, there exists a unique holomorphic solution
$\psi_\lambda^k:\C\rightarrow\C$ of the differential-reflection problem
\begin{equation}\label{deq}
\begin{cases}
T_k f = i\lambda f,\\ f(0) = 1.
\end{cases}
\end{equation}
The map $(z,\lambda)\rightarrow \psi_\lambda^k(z)$ is entire on $\C^2$, and we
have the estimate
\begin{equation}\label{est2}
\left |\psi ^k _\lambda (z)\right | \le e^{|\Im \lambda z|}\qquad(\lambda,z\in \C).
\end{equation}

In view of \eqref{est2} the Dunkl transform $D_k f$ of $f\in L^1 (\R,|w_{k}(x)|dx)$,
where the complex-valued weight function $w_k$ is given by $w_k(x)=  |x|^{2k}$, is
meaningfully defined by
\begin{equation}\label{trafodef}
D_k f (\lambda) = \frac{1}{c_k}\int _\R f(x)\psi ^k _{-\lambda} (x) w_k(x)dx \qquad
( \lambda \in \R),
\end{equation}
where
\[
c_k=\int_\R e^{-\frac{|x|^2}{2}}w_k(x)\,dx\neq 0.
\]
We note that $D_0$ is the Fourier transform on $\R$. From \eqref{est2} we
conclude that $D_k f$ is bounded for such $f$, in fact
\begin{equation}\label{boundedest}
|D_k f(\lambda)|\leq \frac{1}{|c_k|}\int_\R |f(x)|
|w_k(x)|\,dx\quad(\lambda\in\R,\,f\in L^1 (\R,|w_{k}(x)|dx)).
\end{equation}
 The Dunkl transform is a topological isomorphism of $\S(\R)$ onto itself,
the inverse transform $D_k ^{-1}$ being given by
\begin{equation*}
D_k ^{-1} f(x) = \frac{1}{c_k}\int _\R f(\lambda )\psi ^k _{\lambda} (x) w_k
(\lambda)d\lambda=D_k  f(-x) \qquad (f\in \S (\R),\,x \in \R).
\end{equation*}

The operator $T_k$ is anti-symmetric with respect to the weight function $w_k$, i.e.,
\begin{equation}\label{antisym}
\langle T_k f, g\rangle _k = -\langle f, T_k g\rangle _k ,
\end{equation}
for $f\in \S(\R)$ and $g\in C^\infty (\R)$ such that both $g$ and $T_k g$ are of at
most polynomial growth. Here $\langle f, g\rangle_k$ is defined by
 \begin{equation*}
\langle f, g\rangle _k = \int _\R f(x)g(x) w_k(x)dx  ,
\end{equation*}
for functions $f$ and $g$ such that $fg\in L^1 (\R,|w_{k}(x)|dx)$. In particular,
\eqref{antisym} yields the intertwining identity
\[
D_k (T_k f) (\lambda) =i\lambda (D_k f) (\lambda)\quad (f\in \S(\R),\,\lambda \in \R).
\]

Furthermore, for $\lambda,z,s\in\C$ the symmetry properties $\psi ^k
_\lambda(z) =\psi ^k _z(\lambda)$ and $\psi ^k _{s\lambda} (z) =\psi ^k
_\lambda (sz)$ are valid. Using the first of these and \eqref{boundedest}, an
application of Fubini gives
\begin{equation}\label{ex}
\langle D_k f, g\rangle _k = \langle f, D_k g\rangle _k \quad (f,g\in L^1
(\R,|w_{k}(x)|dx)).
\end{equation}
If $k\geq 0$, the Plancherel theorem states that $D_k$ preserves the weighted
two-norm on $L^1 (\R,w_{k}(x)dx)\cap
L^2 (\R,w_{k}(x)dx)$ and extends to a unitary operator on $L ^2(\R,w_k(x)dx)$.

\section{Paley--Wiener theorem for smooth functions}\label{PWfunctions}

The method of proof in this section is inspired by results of Bang \cite{Ba1}.
To be more specific, for $R>0$ let $\H_R(\C)$ denote the space of entire
functions $f$ with the property that, for all $n\in \N\cup\{0\}$, there exists
a constant $C_{n,f} >0$ such that
\begin{equation*}
|f(z)| \le C_{n,f} (1+|z|)^{-n} e^{R |\Im z|}\qquad (z\in \C).
\end{equation*}
Then, if $k\geq 0$ and $f\in\H_R(\C)$, we will establish that (cf.~\cite{Ba1})
\begin{equation}\label{chain}
\sup \{|\lambda| \,:\, \lambda \in \supp D_k f\}\le\liminf_{n \to \infty} \left
\Vert T_k ^n f\right \Vert _\infty ^{1/n} \le \limsup_{n \to \infty} \left
\Vert T_k ^n f\right \Vert _\infty ^{1/n}\le R,
\end{equation}
after which the proof of the Paley--Wiener theorem for smooth functions is a
mere formality.

Starting towards the third of these inequalities, we first use \eqref{intDunkl}
to gain control over repeated Dunkl derivatives.

\begin{lemma}\label{induction}
Let $k\in\C$, $f\in C^\infty (\R)$ and $n\in \N$. Then
\begin{equation*}
T_k ^n f(x) = (T_k ^{n-1} (f')) (x) + k \int_{-1} ^1 t^{n-1} (T_k ^{n-1} ( f'))(tx)
dt \qquad (x\in \R).
\end{equation*}
\end{lemma}
\begin{proof}
For $g\in C^\infty(\R)$ and $m\in\N\cup\{0\}$, let $I_{g,m}\in C^\infty(\R)$ be
defined by
\[
I_{g,m}(x)=\int_{-1}^1 t^m g(tx)\,dt\quad(x\in\R).
\]
Using \eqref{intDunkl}, we then find
\begin{align*}
(T_k I_{g,m})(x)&=\int _{-1}^1 t_1^{m+1}  g' (t_1 x) \,dt _1 +k\int _{-1}^1\int
_{-1}^1 t_1^{m+1} g ' (t_1t_2x) \,dt _1\,dt _2\\
 &= \int _{-1}^1 t_1^{m+1}\left ( g' (t_1 x) +k\int _{-1}^1 g' (t_1t_2 x) \,dt _2
\right )dt _1\\&=\int_{-1}^1 t_1^{m+1} (T_k g)(t_1 x)\, dt_1=I_{T_k g,m+1}(x).
\end{align*}
We conclude that $T_k I_{g,m}=I_{T_k g,m+1}$. Since \eqref{intDunkl} can be written
as $T_k f=f'+k I_{f',0}$ one has
\begin{equation*}
T_k^{n}f=T_k^{n-1}(f'+k I_{f',0})=T_k^{n-1}(f')+k I_{T_k^{n-1}(f'),n-1},
\end{equation*}
which is
the statement in the lemma.
\end{proof}
It follows from Lemma~\ref{induction} that
\begin{equation*}
|T_k ^n f (x)| \le \left (1+\frac {2|k|}{n}\right )  \sup_{y\in [-|x|,|x|]}
|(T_k^{n-1} (f'))(y)| \qquad (x\in \R),
\end{equation*}
and induction then yields the following basic estimate, which is more explicit than
\cite[Prop.~2.1]{CT}.

\begin{corollary}\label{T^n}
Let $k\in\C$, $f\in C^\infty (\R)$ and $n\in \N$. Then
\begin{equation*}
|T_k ^n f (x)| \le \frac {\Gamma(n+1+2|k|)}{n!\,\Gamma(1+2|k|)} \sup_{y\in
[-|x|,|x|]} |f ^{(n)}(y)| \qquad (x\in \R).
\end{equation*}
\end{corollary}
The third inequality in \eqref{chain} can now be settled.

\begin{proposition}\label{exp}
Let $k\in\C,\,R>0$, and suppose $f:\C\to\C$ is an entire function such that
\[
|f(z)|\le C e^{R |\Im z|} \qquad (z\in\C),
\]
for some positive constant $C$. Then, for all $n\in\N$, $T_k^n f$ is bounded on the
real line, and
\begin{equation*}
\limsup _{n \to \infty} \left \Vert T_k ^n f\right \Vert _{\infty} ^{1/n} \le R.
\end{equation*}
\end{proposition}
\begin{proof}
We have, for any $r>0$,
\begin{equation*}
f^{(n)} (z) = \frac {n!}{2\pi i}\oint_{|\zeta - z| = r} \frac {f(\zeta)}{(\zeta
-z)^{n+1}}\, d\zeta \qquad (z\in \C).
\end{equation*}
If $|\zeta - z| =r$, then
\begin{equation*}
|f(\zeta)| \le C e ^{R(|\Im z|+r)},
\end{equation*}
implying
\begin{equation*}
|f^{(n)} (z)| \le C \frac {n!}{r^n}e ^{R(|\Im z|+r)}\qquad (z\in \C).
\end{equation*}
Choosing $r=n/R$, so that $r>0$ if $n\in\N$, we find
\begin{equation*}
|f^{(n)} (z)| \le C \frac {n!e^n}{n^n} R^n e ^{R|\Im z|} \qquad (n\in\N, z\in \C),
\end{equation*}
whence $\Vert f^{(n)}\Vert_\infty\leq Cn! e^n n^{-n}R^n$, for $n\in\N$.
Combining this with Corollary~\ref{T^n} yields
\begin{equation*}
\| T_k^n f\|_\infty \le C
\frac{e^n\,\Gamma(n+1+2|k|)}{n^n\,\Gamma(1+2|k|)}R^n\quad(n\in\N).
\end{equation*}
The result now follows from Stirling's formula.
\end{proof}

As to the first inequality in \eqref{chain}, it is actually easy to prove that
it holds for the norm $\Vert\cdot \Vert_{k,p}$ in $L^p(\R,w_k(x)dx)$ for
arbitrary $1\le p\le \infty$ (not just for $p=\infty$), as is shown by the
following lemma. It should be noted that this result can be generalized---with
different proofs---to complex multiplicities (cf.\ Lemma~\ref{real}) and to
$L^p$-functions for $k\geq 0$ (cf.\ Theorem~\ref{real2}), but we present it
here separately nevertheless, in order to illustrate that for the case $k\geq
0$, the proof of one of the crucial inequalities (as far as the Paley--Wiener
theorem is concerned) is rather elementary and intuitive.

\begin{lemma}\label{easy}
Let $k \ge 0$,  $1\le p\le \infty$ and $f\in \S (\R)$.
Then in the extended positive real numbers, 
\begin{equation}
\liminf_{n \to \infty} \left \Vert T_k ^n f\right \Vert _{k,p} ^{1/n} \ge \sup
\{|\lambda| \,:\, \lambda
\in \supp D_k f\}.
\label{inftyest}
\end{equation}
\end{lemma}

\begin{proof}
Suppose $0\neq\lambda _0\in\supp D _k f$ and let $0<\varepsilon<|\lambda_0|$.
Define
\begin{equation*}
\phi (\lambda) = \overline{D_k f(\lambda)}\qquad (\lambda \in \R).
\end{equation*}

Then, with $q$ denoting the conjugate exponent and using \eqref{ex}, we find
\begin{align*}
\left \Vert T_k ^{2n} f \right \Vert _{k,p}\left \Vert D _k \phi \right \Vert _{k,q}
&\ge
|\langle T_k ^{2n} f, D _k \phi \rangle _k |= |\langle D_k ( T_k ^{2n} f),
\phi\rangle _k |\\
&=\left |\int _{\R} (i\lambda)^{2n} D_k  f(\lambda) \phi (\lambda) w_k
(\lambda)d\lambda \right |\\
&=\int _{\R} \lambda^{2n} |D_k  f(\lambda)|^2 w_{k} (\lambda)d\lambda \\
&\ge (|\lambda _0|-\varepsilon)^{2n} \int _{|\lambda|\ge |\lambda _0| -\varepsilon}
| D_k f(
\lambda)|^2
w_{k} (\lambda)d\lambda.
\end{align*}

With $\psi (\lambda) = \overline{D_k (T_k f)(\lambda) }$
we similarly get
\begin{equation*}
\left \Vert T_k ^{2n+1} f \right \Vert _{k,p}\left \Vert D _k \psi \right \Vert
_{k,q}\ge
(|\lambda _0|-\varepsilon)^{2n+1} \int _{|\lambda|\ge |\lambda _0|
-\varepsilon}|\lambda| | D_k f(\lambda)|^2
w_{k} (\lambda)d\lambda.
\end{equation*}
These two estimates together yield
\begin{equation*}
\liminf _{n \to \infty} \left \Vert T_k ^n f\right \Vert _{k,p} ^{1/n} \ge |\lambda
_0|-\varepsilon,
\end{equation*}
and the lemma follows.
\end{proof}

Now that \eqref{chain} has been established, we come to the Paley--Wiener
theorem for smooth functions. Introducing notation, for $R>0$, we let $C
^\infty _R (\R)$ denote the space of smooth functions on $\R$ with support in
$[-R,R]$. Its counterpart under the Dunkl transform, the space $\H_R(\C)$, was
defined at the beginning of this section.

\begin{theorem}[Paley--Wiener theorem for smooth functions]\label{PWT}
Let $R>0$ and $k \ge 0$. Then the Dunkl transform $D_k$ is a bijection from
$C^\infty _R (\R)$ onto $\H_R(\C)$.
\end{theorem}
\begin{proof}
If $f\in C^\infty _R (\R)$, then it is easy to see that $D_k f\in \H_R(\C)$
\cite[Corollary~4.10]{dJ1}.
Now assume that $f\in \H_R(\C)$. Using Cauchy's integral representation as in the proof
of Proposition~\ref{exp}, we
retrieve the well-known fact that $f\in \S(\R)$. From \eqref{chain} we infer that
$D_k f$ has support in $[-R,R]$. Since $D_k^{-1}f(x)=D_kf(-x)$, for $x\in\R$, the
same is
true for $D_k^{-1}f$, as was to be proved.
\end{proof}

\begin{remark}\label{PWremark}\quad
\begin{enumerate}
\item By holomorphic continuation and continuity, cf.~\cite{dJ2}, one sees that
Theorem~\ref{PWT} also holds in the more general case $\Re k\geq 0$.
Alternatively, one can use Lemma~\ref{real} below instead of Lemma~\ref{easy},
which establishes \eqref{chain} also in the case $\Re k\geq 0$, and then the
above direct proof is again valid. \item We emphasize that the present proof
does not use any explicit formulas for the Dunkl kernel in one dimension,
contrary to the alternative methods of proof in \cite{T} (where Weyl fractional
integral operators are used), \cite{dJ2} (where asymptotic results for Bessel
functions are needed) and \cite{CT} (where various integral operators, Dunkl's
intertwining operator and the Paley--Wiener theorem for the Fourier transform
all play a role). Also, the fact that a contour shifting argument for the
transform is usually not possible (since $w_k$ generically has no entire
extension) is no obstruction. Given this unspecific nature, it is possible that
the present method can be applied to other transforms as well, although the
symmetry of the kernel---as reflected in \eqref{ex}, which was used in the
proof of Lemma~\ref{easy} and which will again be used in the proof of the
alternative Lemma~\ref{real} below---is perhaps necessary. The same suggestion
applies to the results in the remaining sections of this paper.
\end{enumerate}
\end{remark}

\section{Real Paley--Wiener theorem for $L^p$-functions}\label{Lp}

We will now consider the real Paley--Wiener theorem for $L^p$-functions in the
spirit of Bang~\cite{Ba1}. The result is first proved for Schwartz functions in
Theorem~\ref{real1} and subsequently for the general case in Theorem~\ref{real2}.

Let $\Vert\cdot \Vert _{\Re k,p}$ denote the $L^p(\R,|w_k(x)|dx)$-norm, for $1\le
p\le \infty$. Then we have the following generalization of Lemma~\ref{easy} to
complex multiplicities.

\begin{lemma}\label{real}
Let $\Re k \ge 0$,  $1\le p\le \infty$ and $f\in \S (\R)$. Then in the extended
positive real numbers, 
\begin{equation}
\liminf _{n \to \infty} \left \Vert T_k ^n f\right \Vert _{\Re k,p} ^{1/n} \ge \sup
\{|\lambda| \,:\, \lambda
\in \supp D_k f\}.
\label{inftyest2}
\end{equation}
\end{lemma}
\begin{proof}
Let $0\ne \lambda _0 \in\supp D _k f$ and choose $\epsilon>0$ such that
$0<2\varepsilon <|\lambda_0|$. Also choose $\phi\in C_c ^\infty(\R)$ such that
$\supp \phi \subset
[\lambda _0 -\varepsilon,\lambda _0 +\varepsilon]$, and
$\langle D_k f, \phi \rangle_k  \ne 0$.
Define $\phi _{n} (\lambda) = \lambda  ^{-n} \phi(\lambda )$ and $P_n
(x) = x^n$ for $n\in \N\cup\{0\}$.
Then
\begin{equation*}
(1+ P_{N}(x)) (D_k \phi _{n}) (x) =\frac{1}{c_k}
\int _{\lambda _0 -\varepsilon } ^{\lambda _0 +\varepsilon } \left (1+ (iT _k)^
N \right )
( \lambda  ^{-n} \phi (\lambda )) \psi _x^k  (\lambda) w_k (\lambda) d \lambda
\quad(N\in \N\cup\{0\}).
\end{equation*}
We fix $N$ such that $N$ is even and $N>2\,\Re k +1$.

Corollary~\ref{T^n} and the binomial formula imply that
\begin{equation*}
\left | \left (1 +(iT_k) ^N\right )
( \lambda  ^{-n} \phi (\lambda ))\right |\le C_1 n^{N}
({|\lambda _0| -\varepsilon} )^{-n}\qquad (n\in \N\cup\{0\},\, \lambda  \in \R),
\end{equation*}
where
$C_1$ is a positive constant. This yields the estimates
\begin{align*}
\Vert  D_k  \phi _{n} \Vert _{\Re k,q} &\le
\Vert (1+ P_{N})^{-1} \Vert _{\Re k,q}\Vert  (1+ P_{N}) D_k \phi _{n}  \Vert _\infty \\
&\le  \frac{2\varepsilon}{|c_k|} C_1
n^{N}\left (|\lambda _0| -\varepsilon \right )^{-n}
\Vert (1+ P_{N})^{-1} \Vert _{\Re k,q}\\&\le C_2
n^{N}\left (|\lambda _0| -\varepsilon \right )^{-n}
\end{align*}
for all $n > N$, where $C_2$ is a positive constant and $q$ is the conjugate exponent.

Using \eqref{ex}, the identity
$D_k (P _n \phi _n) = (-i)^n T_k ^n D_k \phi _n$ and H\"older's inequality, we
therefore get
\begin{align*}
 |\langle D_k f, \phi \rangle _k| &=
 |\langle D_k f, P_{n}  \phi _{n} \rangle_k | = |\langle  f, D_k (P_{n}  \phi _{n})
 \rangle_k |
 =|\langle  f, T _k ^n (D_k \phi _{n})\rangle_k |\\&
= \left| \left\langle T _k ^n f, D_k \phi _{n}\right\rangle_ k \right|\le
\left \Vert T _k ^n f\right \Vert _{\Re k,p}
\left \Vert  D_k \phi _{n} \right \Vert _{\Re k,q}\\&
\le C_2 n^{N}\left (|\lambda _0| -\varepsilon \right )^{-n}
\left \Vert T _k ^n f\right \Vert _{\Re k,p},
\end{align*}
whence
\begin{equation*}
\liminf _{n \to \infty}\left \Vert T _k ^nf\right \Vert _{\Re k,p} ^{1/n}
\ge \liminf _{n \to \infty}\,
(C_2 n^{N}) ^{-1/n}\left (|\lambda _0| -\varepsilon \right )|\langle D_k f, \phi
\rangle_k|^{1/n}
=\left (|\lambda _0| -\varepsilon \right ),
\end{equation*}
establishing the lemma.
\end{proof}

Using the Paley--Wiener theorem, we can extend the inequality in
Proposition~{\ref{exp}} to the norms
$\Vert\cdot \Vert _{\Re k,p}$, $1\le p\le \infty$, for $\Re k\ge 0$,
and we thus have the following theorem.

\begin{theorem}[real Paley--Wiener theorem for Schwartz functions]\label{real1}
Let $\Re k \ge 0$, $1\le p\le\infty$ and $f\in \S (\R)$. Then in the extended
positive real numbers, 
\begin{equation*}
\lim _{n \to \infty} \left \Vert T_k ^n f\right \Vert _{\Re k,p} ^{1/n} = \sup
\{|\lambda| \,:\, \lambda
\in \supp D_k f\} .
\end{equation*}
\end{theorem}
\begin{proof}
In view of Lemma~\ref{real} it only remains to be shown that
\begin{equation*}
\limsup _{n \to \infty} \left \Vert T_k ^n f\right \Vert _{\Re k,p} ^{1/n} \le R,
\end{equation*}
if $f\in \S(\R)$ is such that $\supp D_k f \subset [-R,R]$ for some finite $R>0$.
Using the inversion formula and the intertwining properties of the transform  we
have
\begin{equation}\label{eq1real1}
x^N T_k^n f (x) =  \frac{i^{N+n}}{c_k}  \int_{-R}^R T_k ^N (P_n D_k f ) (\lambda)
\psi^k _\lambda (x) w_k(\lambda)\,d\lambda \quad (n,N\in\N\cup\{0\}),
\end{equation}
where again $P_n(x)=x^n$. Now Corollary~\ref{T^n} and the binomial formula imply that
\begin{equation*}
\|T_k ^N (P_n D_k f )\|_\infty \le C_1 n ^N R^n,
\end{equation*}
where $C_1$ is a constant depending on $f$ and $N$. Therefore \eqref{eq1real1}
yields that
\begin{equation}\label{eq2real1}
\Vert(1+P_N)T_k^n f \Vert_\infty\le C_2 n ^N R^{n+1},
\end{equation}
where $C_2$ is again a constant depending on $f$ and $N$. We fix $N$ such that $N$
is even and $N>2\, \Re k +1$. Then the observation
\[
\Vert T_k^n f\Vert_{\Re k,p}\leq \Vert (1+P_N)^{-1}\Vert_{\Re k,p}\Vert
(1+P_N)T_k^nf\Vert_\infty
\]
and \eqref{eq2real1} establish the result.
\end{proof}

\begin{remark} Theorem~\ref{real1} is new for complex $k$. For real $k$, the result
can be found in \cite{CT}, where it is proved using the Plancherel theorem for the
Dunkl transform, the Riesz--Thorin convexity theorem and the theory of Sobolev
spaces for Dunkl operators.
\end{remark}

For $k\geq 0$, we will now generalize Theorem~\ref{real1} to the $L^p$-case in
Theorem~\ref{real2},
using the structure of $\S(\R)$ as an associative algebra under the Dunkl
convolution $*_k$. We refer to \cite{S} for
details on this subject.

Let $k\ge 0$, and define a distributional Dunkl transform
$D_k^d f$ of $f\in L ^p(\R,w_{k}(x)dx)$ by transposition
\begin{equation*}
\langle D_k^d f,\phi\rangle =\langle  f, D_k \phi\rangle_k  \qquad (\phi \in \S(\R)).
\end{equation*}
Clearly $D_k^d$ is injective on
$L ^p(\R,w_{k}(x)dx)$. Furthermore, from
\begin{equation*}
|\langle D_k ^d f,\phi\rangle| \le \|  f\|_p\| D_k \phi\|_\infty \qquad (\phi \in
\S(\R)),
\end{equation*}
we see that $D_k ^d f$ is a tempered distribution.
Note also that for $f\in L^p (\R,w_{k}(x)dx)$ and $1\le p\le 2$ we have
\begin{equation*}
\langle D_k^d  f,\phi\rangle =\langle D_k  f,\phi\rangle _k  \qquad (\phi \in \S(\R)).
\end{equation*}
by (\ref{ex}) and density of $\S(\R)$ in $L^p (\R,w_{k}(x)dx)$. Thus, for $f\in L^p
(\R,w_{k}(x)dx)$ and $1\le p\le 2$,  $D_k^d f$ as defined above corresponds to the
distribution $(D_k f) w_k$, implying that in this case $\supp D_k ^d f = \supp D_k f$.

\begin{theorem}[real Paley--Wiener theorem for $L^p$-functions]\label{real2}
Let $k\ge 0$, $1\le p\le \infty$ and $f\in C^\infty(\R)$ be such that
$T_k ^n f \in L^p(\R,w_{k}(x)dx)$, for all $n\in \N\cup\{0\}$. Then in the extended
positive real numbers,
\begin{equation*}
\lim _{n \to \infty} \left \Vert T_k ^n f\right \Vert _{k,p} ^{1/n} = \sup
\{|\lambda| \,:\, \lambda
\in \supp D_k^d f\} .
\end{equation*}
If in addition $f\in L^s(\R,w_{k}(x)dx)$, for some $1\le s\le 2$, then the distribution
$D_k^d f$ corresponds to the function $(D_k f) w_k$, and the support of
$D_k^d f$ in the right hand side is equal to the support of $D_k f$ as a
distribution.
\end{theorem}
\begin{proof}
First we note that \eqref{inftyest2} also holds for $f$ as above: in the proof of
Lemma~{\ref{real}} we just have to change $\langle D_k f,\cdot\rangle _k$ into
$\langle D_k ^df,\cdot\rangle$.
Therefore, it only remains to be shown that
\begin{equation}\label{real2eq1}
\limsup _{n \to \infty} \left \Vert T_k ^n f\right \Vert _{k,p} ^{1/n} \le R,
\end{equation}
if $f$ is as in the theorem and such that $\supp D_k^d f \subset [-R,R]$ for some
finite $R>0$.
To this end, choose $\varepsilon>0$, and fix a function $\phi _{\varepsilon}\in
\S(\R)$ such that
$D_k^{-1}\phi _{\varepsilon}=1$ on $[-R,R]$ and
$D_k^{-1} \phi_{\varepsilon}=0$ outside $[-R-\varepsilon,R+\varepsilon]$. We
have from \cite[Proposition 3]{S} that $D_k^{-1}(\phi*_k\psi)=(D_k^{-1}\phi )\,
(D_k^{-1}\psi)$, for all $\phi,\psi\in\S(\R)$. With $P_n(x)=x^n$,  we thus find for
arbitrary $\phi\in\S(\R)$ that
\begin{align*}
\langle f,\phi*_kT_k^n\phi_\varepsilon\rangle_k&=\langle D_k^d f,D_k^{-1}(\phi*_k
T_k^n\phi_\varepsilon)\rangle=(-i)^n\langle D_k^d f,P_n (D_k^{-1}\phi)\,
(D_k^{-1}\phi_\varepsilon)\rangle\\&=(-i)^n\langle D_k^d f,P_n D_k^{-1}\phi
\rangle=\langle f, T_k^n \phi\rangle_k.
\end{align*}
Furthermore, from loc.cit., one knows that $\|  \phi *_k \psi\|_{k,q} \le 4 \|\phi
\|_{k,q} \|\psi\|_{k,1}$,
for all $\phi,\psi \in \S(\R)$, where $q$ is the conjugate exponent of $p$. If we
combine these two results with \eqref{antisym} and H\"older's reverse inequality,
we infer that
\begin{equation*}
\|T_k ^n f\|_{k,p} = \sup _{\phi}|\langle T_k ^n f,\phi\rangle_k |= \sup
_{\phi}|\langle f,T_k^n \phi\rangle_k |
= \sup _{\phi}|\langle  f,\phi * _k T_k ^n \phi _\varepsilon\rangle_k | \le 4
\| f\|_{k,p} \| T_k ^n \phi _\varepsilon\|_{k,1},
\end{equation*}
where the supremum is over all functions $\phi \in \S(\R)$ with
$\|\phi \|_{k,q} = 1$. From Theorem~\ref{real1} we therefore conclude that
\[
\limsup _{n \to \infty} \left \Vert T_k ^n f\right \Vert _{k,p} ^{1/n} \le
R+\varepsilon,
\]
proving \eqref{real2eq1}.

The statement on supports was established in the discussion preceding the theorem.
\end{proof}

\begin{remark}
For even functions, when the Dunkl transform reduces to the Hankel transform, the
previous result can already be found in \cite{An2, An3}. Also using Dunkl convolution,
and closely following the approach in \cite{An1, An2, An3}, Theorem~\ref{real2} has
previously been established in \cite{COT}. Our proof is considerably shorter than
the proof in loc.cit.
\end{remark}

\section{Paley--Wiener theorem for $L^2$-functions (for $k\ge 0$)}\label{L2}

In this section we assume that $k\ge 0$. For $R>0$ and $1\leq p\leq\infty$, we
define $L_{R} ^p(\R,w_{k}(x)dx)$ to be the subspace of $L ^p(\R, w_{k}(x)dx)$
consisting of those functions with distributional support in $[-R,R]$, and we
let $\H^{p,k}_{R} (\C)$ denote the space of entire functions $f:\C\to\C$ which
belong to $L^{p}(\R,w_{k}(x)dx)$ when restricted to the real line and which are
such that
\[
|f(z)|\le C_f e^{R |\Im z|} \qquad (z\in\C),
\]
for some positive constant $C_f$.

\begin{theorem}[Paley--Wiener theorem for $L^2$-functions]\label{PWL2}
Let $R>0$ and $k \ge 0$. Then the Dunkl transform $D_k$ is a bijection from
$L_{R} ^2(\R,w_{k}(x)dx)$ onto $\H_{R}^{2,k}(\C)$.
\end{theorem}
\begin{proof}
Let $f\in L_{R} ^2(\R,w_{k}(x)dx)$. Then $f\in L^1(\R,w_{k}(x)dx)$, and
(\ref{est2}) and (\ref{trafodef}) together with the Plancherel theorem imply
that $D_k f\in\H^{2,k}_{R} (\C)$. Conversely, let $f\in \H_{R}^{2,k}(\C)$. By
the Plancherel theorem one has $D_k^{-1} f\in L ^2(\R,w_{k}(x)dx)$. In
addition, Proposition~{\ref{exp}} and Theorem~{\ref{real2}} (with $p=\infty$)
show that $\supp D_k f\subset[-R,R]$. The same is then true for $D_k^{-1}f$,
and the result follows.
\end{proof}

\begin{remark} \quad
\begin{enumerate}
\item  If $f\in\H^{p,k}_{R} (\C)$ and $1\leq p\leq\infty$, then using
Proposition~{\ref{exp}} and Theorem~{\ref{real2}} as in the proof of
Theorem~\ref{PWL2}, one sees that $D_k ^d f$ has support in $[-R,R]$. In
particular, for $1\le p\le 2$ with conjugate exponent $q$, we conclude that
$D_k$ maps $\H^{p,k}_{R} (\C)$ into $L_{R} ^q(\R,w_{k}(x)dx)$. For even
functions, when the Dunkl transform can be identified with the Hankel
transform, the latter result can be found in \cite{G}.

\item As an ingredient for the discussion of the relation with the literature
on the Fourier transform, let us make the preliminary observation that, for
$R>0$ and $1\leq p\leq\infty$, an entire function $f$ is in $H^{p,0}_R(\C)$ if,
and only if, its restriction to the real line is in $L^p(\R,dx)$ and moreover
\[
|f(z)|\le \widetilde C_f e^{R |z|} \qquad (z\in\C),
\]
for some positive constant $\widetilde C_f$ \cite[Theorems~6.2.4 and 6.7.1]{Boas}. This
being said, for $p=1$ the specialization of the first part of this remark to
$k=0$ therefore proves part of the statement in \cite[Theorem~6.8.11]{Boas}, and for
$1<p<2$ this specialization proves the first statement of
\cite[Theorem~6.8.13]{Boas}.

\item The aforementioned result about entire functions shows that the
specialization to $k=0$ of Theorem~\ref{PWL2} is equivalent to the original
Paley--Wiener theorem for the Fourier transform (see \cite{PW} or \cite[Theorem
19.3]{R}). The present proof seems to be more in terms of general principles
than other proofs seen in the literature.
\end{enumerate}
\end{remark}

\end{document}